\documentclass[12pt,leqno]{amsart}
\usepackage{amssymb}
\usepackage{amsmath,amssymb,color}
\numberwithin{equation}{section}

\newcommand{\MLML}{(t-s)^{\alpha-1}\MLAA (-\la_n(t-s)^{\alpha})}

\newcommand{\la}{\lambda}
\newcommand{\va}{\varphi}
\newcommand{\ppp}{\partial}

\newcommand{\pppa}{\partial_t^{\alpha}}

\newcommand{\C}{\mathbb{C}} 
\newcommand{\N}{\mathbb{N}}

\newcommand{\OOO}{\Omega}
\newcommand{\MLO}{E_{\alpha,1}}

\newcommand{\MLAA}{E_{\alpha,\alpha}}
\newcommand{\sumn}{\sum_{n=1}^{\infty}}

\newcommand{\CCO}{{_{0}C^1[0,T]}}

\newcommand{\hhalf}{\frac{1}{2}}

\allowdisplaybreaks
\title
[]
{
Uniqueness for an inverse coefficient problem for a one-dimensional 
time-fractional diffusion equation with non-zero boundary conditions
}

\author{
$^1\;$William Rundell
$^{2,3,4}\;$Masahiro Yamamoto }

\thanks{
$^1\;$Department of Mathematics, Texas A\&M University, College Station,
Texas 77843 USA
e-mail: {\tt rundell@math.tamu.edu}
\\
$^2\;$Graduate School of Mathematical Sciences, The University
of Tokyo, Komaba, Meguro, Tokyo 153-8914, Japan \\
$^3\;$Honorary Member of Academy of Romanian Scientists, 
Splaiul Independentei Street, no 54,
050094 Bucharest Romania \\
$^4\;$Peoples' Friendship University of Russia 
(RUDN University) 6 Miklukho-Maklaya St, Moscow, 117198, Russian Federation
e-mail: {\tt myama@ms.u-tokyo.ac.jp}
}

\date{}
\begin{document}
\maketitle

\begin{abstract}
We consider initial boundary value problems for one-dimensional diffusion 
equation with time-fractional derivative of order $\alpha \in (0,1)$ which 
are subject to non-zero Neumann boundary conditions.
We prove the uniqueness for an inverse coefficient problem of 
determining a spatially varying potential and the order of the 
time-fractional derivative by Dirichlet data at one end point of the spatial
interval.  The imposed Neumann conditions are required to be within the 
correct Sobolev space of order $\alpha$.  Our proof is based on 
a representation formula of solution to an initial boundary value problem 
with non-zero boundary data.  Moreover, we apply such a formula and prove
the uniqueness in the determination of boundary value at another end point 
by Cauchy data at one end point. 
\end{abstract}

\section{Introduction}

We consider the following initial boundary value problem for 
a one-dimensional time-fractional diffusion equation:
\begin{equation}\label{eqn:basic_ode}
\left\{ \begin{array}{rl}
& d_t^{\alpha} u(x,t) = \ppp_x^2u(x,t) + p(x)u(x,t), \quad 
0<x<1, \, 0<t<T, \\
& \ppp_xu(0,t) = 0, \quad \ppp_xu(1,t) = g(t), \quad 0<t<T,\\
& u(x,0) = 0, \quad 0<x<1.
\end{array}\right.  %\eqno{(1.1)}
\end{equation}
Here and henceforth let $\ppp_x = \frac{\ppp}{\ppp x}$, 
$\ppp_x^2 = \frac{\ppp^2}{\ppp x^2}$,
and we define for absolutely continuous $g$ on $[0,T]$
$$
d_t^{\alpha} g(t) = \frac{1}{\Gamma(1-\alpha)}\int^t_0
(t-s)^{-\alpha}\frac{dg}{ds}(s)\, ds, \quad 0<t<T, 
$$
that is, the fractional derivative of order $\alpha$,
$0<\alpha < 1$, 
and of Caputo type (see, for example, Podlubny \cite{Po}).
The first equation in \eqref{eqn:basic_ode}
%(1.1)
is a time-fractional diffusion equation of subdiffusion type
modelling, for example, anomalous diffusion in heterogeneous 
media. For some applications see, for example, Metzler and Klafter \cite{MK}.

In this article we are concerned with the question of uniqueness
for the inverse problem:\\
{\it Let $g = g(t)$ be given for $0<t<T$.  
Given data $u(0,t)$ for $0<t<T$ or $u(1,t)$ for $0<t<T$,  does this
uniquely determine 
$\; 
\alpha \in (0,1) \quad \mbox{and} \quad p(x), \, 0<x<1?
\;$
}
\\

In place of \eqref{eqn:basic_ode} we can also consider 
\begin{equation}\label{eqn:basic_ode_2}
\left\{ \begin{array}{rl}
& d_t^{\alpha} u(x,t) = \ppp_x^2u(x,t) + p(x)u(x,t), \quad 
0<x<1, \, 0<t<T, \\
& \ppp_xu(0,t) = \ppp_xu(1,t) = 0, \quad 0<t<T, \\
& u(x,0) = a(x), \quad 0<x<1.
\end{array}\right.
\end{equation}
Uniqueness for this type of inverse problem for \eqref{eqn:basic_ode_2} with 
$\alpha=1$, that is, for the initial boundary value problem for the 
heat equation, was considered by, for example, Murayama \cite{Mu}, 
Suzuki and Murayama \cite{SM}.
For the case with 
$0<\alpha<1$, we refer to Cheng, Nakagawa, Yamamoto and Yamazaki \cite{CNYY},
Li, Zhang, Jia and Yamamoto \cite{LZJY}.  Also see Jin and Rundell 
\cite{JR}, Jing and Yamamoto \cite{JY}, and survey chapters 
Li, Liu and Yamamoto \cite{LiLiuYa},
Li and Yamamoto \cite{LiYa}, Liu, Li and Yamamoto \cite{LiuLiYa}.
Both for the cases of $\alpha=1$ and $0<\alpha<1$, the uniqueness for
\eqref{eqn:basic_ode_2}
requires a quite strong condition to be imposed for the initial value $a(x)$.

On the other hand, for the inverse problem for \eqref{eqn:basic_ode}
with a zero initial value but $g \not\equiv 0$, 
we refer to Pierce \cite{Pi} who proved the uniqueness 
for $\alpha=1$ with the quite mild assumption $g\not \equiv 0$.

For fixed $\alpha \in (1,2)$, Wei and Yan \cite{WY} established 
the uniqueness in determining $p(x)$ with $g \in C^2[0,T]$ 
imposing additional conditions.

For the inverse problem for \eqref{eqn:basic_ode} with $0<\alpha<1$,
see Rundell and Yamamoto \cite{RY}.  The purpose of this article is 
to complete \cite{RY} within a weaker class of solutions in 
suitable Sobolev space in time.  For the case of $1<\alpha<2$, we can 
argue in a similar manner but we concentrate on the case $0<\alpha<1$.
\\

For the mathematical formulations, we need to introduce function spaces and 
relevant operators;
all functions considered are assumed to be real-valued.
Let $L^2(0,1)$ be a usual Lebesgue space and let $\langle \cdot,\cdot\rangle$ 
and $\Vert \cdot\Vert$ denote the scalar product and the norm respectively
in $L^2(0,1)$, and let $\langle \cdot,\cdot \rangle_X$ be the scalar product 
in other Hilbert spaces $X$ when we so specify.

We define the fractional Sobolev space $H^{\alpha}(0,T)$ 
on the interval $(0,T)$ (see e.g., \cite{Ad}, Chapter VII) with the norm 
in $H^{\alpha}(0,T)$: 
$$
\Vert u\Vert_{H^{\alpha}(0,T)}
:= \left( \Vert u\Vert^2_{L^2(0,T)}
+ \int^T_0 \int^T_0 
\frac{\vert u(t)-u(s)\vert^2}{\vert t-s\vert^{1+2\alpha}}\,dtds
\right)^{\!\frac{1}{2}}. 
$$
We further define the Banach spaces
$$
H_{\alpha}(0,T) := 
\left\{ \begin{array}{rl}
& \{ u \in H^{\alpha}(0,T); \thinspace
u(0) = 0\}, \quad  \frac{1}{2} < \alpha < 1, \\
&\left\{v \in H^{\frac{1}{2}}(0,T); \thinspace \int^T_0 
\frac{\vert v(t)\vert^2}
{t}\, dt < \infty\right\}, \quad \alpha=\frac{1}{2}, \\
& H^{\alpha}(0,T), \quad  0 < \alpha < \frac{1}{2}
\end{array}\right.
$$ 
with the following norm:
$$
\Vert v\Vert_{H_{\alpha}(0,T)}=
\left\{ \begin{array}{rl}
\Vert v\Vert_{H^{\alpha}(0,T)}, \quad  &0<\alpha< 1, \thinspace \alpha
\ne \frac{1}{2}, \\
\left( \Vert v\Vert^2_{H^{\hhalf}(0,T)} 
+ \int^T_0 \frac{\vert v(t)\vert^2}{t}\, dt \right)^{\frac{1}{2}}, \quad
& \alpha = \hhalf.  
\end{array}\right.              
$$
We define the Abel (Riemann-Liouville) fractional integral operator 
\begin{equation}
%\label{eqn:Abel_int}
J^{\alpha}g(t) = \frac{1}{\Gamma(\alpha)}\int^t_0
(t-s)^{\alpha-1} g(s)\, ds, \quad 0<t<T,\quad 0<\alpha<1.
\nonumber
\end{equation}
Henceforth by $x \sim y$, we mean that there exists a constant
$C>0$ such that $C^{-1}y \le x \le Cy$ for all quantities $x, y$ 
under consideration.
 
In Gorenflo, Luchko and Yamamoto \cite{GLY}, 
Kubica, Ryszewska and Yamamoto \cite{KRY} (Theorem 2.1), it is proved that 
$J^{\alpha}$ is an isomorphism between $L^2(0,T)$ and $H_{\alpha}(0,T)$.
We define
$$
\pppa g = (J^{\alpha})^{-1}g \quad \mbox{for $g \in H_{\alpha}(0,T)
= J^{\alpha}L^2(0,T)$}.
$$
Then also by Theorem 2.5 in \cite{KRY}, we see
\begin{equation}\label{eqn:g_conditions} %\eqno{(1.3)}
\left\{ \begin{array}{rl}
& \Vert \pppa g\Vert_{L^2(0,T)} \sim \Vert g\Vert_{H_{\alpha}(0,T)},
\quad g\in H_{\alpha}(0,T), \\
& \pppa g = d_t^{\alpha}g \quad \mbox{if 
$g \in W^{1,1}(0,T)$ satisfies $g(0) = 0$
and $t^{\alpha-1}\frac{dg}{dt} \in L^{\infty}(0,T)$}.
\end{array}\right.
\end{equation}
In other words,
$\pppa$ is an extension of the Caputo derivative $d_t^{\alpha}$
to $H_{\alpha}(0,T)$.

Thus throughout this article, in place of \eqref{eqn:basic_ode} we consider 
\begin{equation}\label{eqn:basic_ode_3} %\eqno{(1.4)}
\left\{ \begin{array}{rl}
& \pppa u(x,t) = \ppp_x^2u(x,t) + p(x)u(x,t), \quad 
0<x<1, \, 0<t<T, \\
& \ppp_xu(0,t) = 0, \quad \ppp_xu(1,t) = g(t), \quad 0<t<T,\\
& u \in H_{\alpha}(0,T;L^2(0,1)).
\end{array}\right.
\end{equation}
We assume
\begin{equation}\label{eqn:p_q_negative}  %\eqno{(1.5)}
p,\, q \le 0 \quad \mbox{on $[0,1]$}, \quad 
p, q \in C[0,1]. 
\end{equation}
Then we can prove
\\
{\bf Proposition 1.}\\
{\it
Let $g \in H_{\alpha}(0,T)$ and let $0<\alpha<1$.
Then there exists a unique solution $u_{p,\alpha}=u_{p,\alpha}(x,t)
\in H_{\alpha}(0,T;L^2(0,1)) \cap L^2(0,T;H^2(0,1))$ solving
\eqref{eqn:basic_ode_3}.
}
\\

In \eqref{eqn:basic_ode_3}, we interpret $u(x,\cdot) \in H_{\alpha}(0,T)$ as 
an initial condition: if $\alpha>\hhalf$, then the Sobolev embedding
yields $H_{\alpha}(0,T;L^2(0,1)) \subset H^{\alpha}(0,T;L^2(0,1)) \subset
C([0,T];L^2(0,1))$ and so this means that $u$ satisfies the initial condition 
in a usual sense.  However for $\alpha<\hhalf$, the time regularity 
does not admit such a usual initial condition and alternatively 
the third equation in \eqref{eqn:basic_ode_3} is required. 
For the class of solutions with the $H_{\alpha}$-regularity in 
$t$, it is sufficient to assume the same regularity in $t$ for boundary 
data $g(t)$, 
that is, $g \in H_{\alpha}(0,T)$.  Moreover for $\alpha>\hhalf$, 
the condition means that $g(0) = 0$, which is a natural compatibility
condition at $x=0$ and $t=0$.  We emphasize that since the order of time 
derivative appearing in the equation is up to $\alpha<1$, it is natural to
work within "$\alpha$-time differentiability", and not in the $C^1$ nor 
$H^1$-class.

For the initial boundary value problems with the zero boundary values,
we refer to Gorenflo, Luchko and Yamamoto \cite{GLY},
Kian and Yamamoto \cite{KiYa}, Kubica, Ryszewska and Yamamoto \cite{KRY},
Kubica and Yamamoto \cite{KY}, Luchko \cite{Lu}, Sakamoto and Yamamoto
\cite{SY}.  On the other hand, for initial boundary value problems with 
non-zero boundary data, there are not many works and we refer only to 
Yamamoto \cite{Ya} in the case of less regular boundary data, and 
one can consult the references therein.  On the other hand, 
the proof of Proposition~1 can be done directly,
thanks to the one-dimensionality, and see Section 2.
\\

Now we are ready to state the main result of this article.
\\
{\bf Theorem 2.}\quad
{\it We assume \eqref{eqn:p_q_negative} and $0< \alpha, \beta < 1$,
\begin{equation}\label{1.6}
g \in H_{\max \{\alpha,\beta\}}(0,T), \quad 
g\not\equiv 0 \quad \mbox{in $(0,T)$}.
\end{equation}
Then either $u_{p,\alpha}(0,t) = u_{q,\beta}(0,t)$ for $0<t<T$ or
$u_{p,\alpha}(1,t) = u_{q,\beta}(1,t)$ for $0<t<T$, yields
%}
%\vskip-1.4\baselineskip
$$
\alpha =\beta, \quad p(x) = q(x), \quad 0<x<1.
$$
}
\\

By the regularity shown in Proposition~1 and the trace theorem, we notice that 
the data $u_{p,\alpha}(0,t)$, etc. can make sense in $L^2(0,T)$.
We stress that the condition $g\not\equiv 0$ in \eqref{1.6}
for the boundary input is quite generous.

The article is composed of four sections.  In Section 2, we prove 
Proposition~1 and a key representation formula of the solution 
$u_{p,\alpha}$ to (1.4). 
Section 3 is devoted to the proof of Theorem 2 
on the basis of the representation formula in Section 2.
In Section 4, we provide one application of the representation formula
to prove the uniqueness in determining a boundary value
at $x=1$ by by Cauchy data at $x=0$. 
\\

\section{Proof of Proposition 1 and a representation formula}

{\bf 2.1. Proof of Proposition~1.}

Recalling \eqref{eqn:p_q_negative}
we define an operator $A_p$ in $L^2(0,1)$ by 
$$
\left\{ \begin{array}{rl}
&A_pw(x) = - \frac{d^2w}{dx^2}(x) - p(x)w(x), \quad 0<x<1,\\
& \mathcal{D}(A_p) = \left\{ w \in H^2(0,1);\,
\frac{dw}{dx}(0) = \frac{dw}{dx}(1) = 0 \right\}.
\end{array}\right.
$$
Then $A_p$ possesses eigenvalues
$
0 < \la_1 < \la_2 < \cdots.
\;$
Let $\va_n$, $n\in \N$ be the associated unique eigenfunction
for $\la_n$: $\va_n
\in \mathcal{D}(A_p)$ satisfies
$A_p\va_n = \la_n\va_n$ in $(0,1)$ and we make the normalisation 
$\va_n(1) = 1$.
Moreover, it is known that $\langle \va_n, \va_m \rangle :=
\int_{\OOO} \va_n(x)\va_m(x) dx = 0$ for $n\ne m$ and we
set the associated norming constants as
$$
\rho_n := \Vert \va_n\Vert^2, \quad n\in \N.
$$

We define
\begin{equation}\label{eqn:v-f}  %\eqno{(2.1)}
\left\{ \begin{array}{rl}
&v(x,t) = u_{p,\alpha}(x,t) - \frac{x^2}{2}g(t), \\
& f(x,t) = -\frac{x^2}{2}\pppa g(t) + g(t) + \frac{x^2}{2}p(x)g(t),
\quad 0<x<1,\, 0<t<t.   
\end{array}\right.
\end{equation}
Then (1.4) is equivalent to
\begin{equation}\label{eqn:basic_ode_4} %\eqno{(2.2)}
\left\{ \begin{array}{rl}
& \pppa v(x,t) = \ppp_x^2v(x,t) + p(x)v(x,t) + f(x,t), \quad 
0<x<1, \, 0<t<T, \\
& \ppp_xv(0,t) = \ppp_xv(1,t) = 0, \quad 0<t<T,\\
& v \in H_{\alpha}(0,T;L^2(0,1)).
\end{array}\right.
\end{equation}
Since $g \in H_{\alpha}(0,T)$, we see that 
$v \in H_{\alpha}(0,T;L^2(0,1))$ if and only if
$u \in H_{\alpha}(0,T;L^2(0,1))$.
\\

From $g \in H_{\alpha}(0,T)$ and $p\in C[0,1]$, it follows that 
$f\in L^2(0,T;L^2(0,1))$.  Thus it is sufficient to prove the 
unique existence of solution $v \in 
H_{\alpha}(0,T;L^2(0,1)) \cap L^2(0,T;H^2(0,1))$ to (2.2).
This follows from \cite{KRY}, \cite{SY} for example.
We note that in \cite{KRY}, \cite{SY}, the zero Dirichlet boundary 
condition is considered and the case of he zero Neumann 
boundary condition can be treated in the same way.
Thus the proof of Proposition~1 is complete.
$\square$
\\

\noindent
{\bf 2.2. The representation formula.}

For $\gamma_1, \gamma_2 > 0$, we define the two parameter
Mittag-Leffler function:
$$
E_{\gamma_1, \gamma_2}(z) = \sum_{k=0}^{\infty}
\frac{z^k}{\Gamma(\gamma_1k + \gamma_2)}, \quad z \in \C.
$$
This is an entire function of order 1 in $z \in \C$ (e.g., Gorenflo, Kilbas,
Mainardi and Rogosin \cite{GKMR}, Podlubny \cite{Po}).
Then
\\
{\bf Proposition 3 (representation formula).}\\
{\it
Let $0<\alpha<1$, $p$ satisfy (1.5) and $g \in H_{\alpha}(0,T)$.
Then
\begin{equation}\label{eqn:u_p,alpha-representation} %\eqno{(2.3)}
u_{p,\alpha}(x,t) = \sumn \frac{1}{\rho_n}\left(\int^t_0
(t-s)^{\alpha-1}\MLAA (-\la_n(t-s)^{\alpha}) g(s)\,ds\right) \, \va_n(x)
\end{equation}
in $H_{\alpha}(0,T;L^2(0,1)) \cap L^2(0,T;H^2(0,1))$.
}
\\
{\bf Proof of Proposition 3.}

By \cite{SY} for example, we have the representation
\begin{equation}\label{eqn:v-representation} %\eqno{(2.4)}
v(x,t) = \sumn \frac{1}{\rho_n}\left(\int^t_0
\MLML \langle f(\cdot,s), \va_n\rangle\, ds\right) \, \va_n(x) 
\end{equation}
in $H_{\alpha}(0,T;L^2(0,1)) \cap L^2(0,T;H^2(0,1))$.
Here we note equations \eqref{eqn:v-f} and (2.2). 

We set
$\;
\CCO := \{ h\in C^1[0,T];\, h(0) = 0\}.
$
\\
First we prove \eqref{eqn:u_p,alpha-representation} for $g \in \CCO$.
We have to calculate the right-hand side of \eqref{eqn:v-representation}.
\begin{equation}\label{eqn:rhs_v-representation} %\eqno{(2.5)}
\begin{aligned}
\int^t_0 &\MLML \langle f(\cdot,s), \va_n\rangle\, ds \\
& = - \int^t_0 \MLML \ppp_s^{\alpha}g(s)\, ds\,
\bigl\langle\frac{x^2}{2},\va_n\bigr\rangle \\
&\quad + \int^t_0 \MLML g(s)\, ds\,
\bigl\langle1+\frac{x^2}{2}p,\, \va_n\bigr\rangle.  \\
\end{aligned}
\end{equation}
We set 
$$
S:= \int^t_0 \MLML \ppp_s^{\alpha}g(s)\, ds.
$$
For $g \in \CCO$, by \eqref{eqn:g_conditions} we see 
that $\ppp_s^{\alpha}g$ coincides with $d_s^{\alpha}g$:
$$
\ppp_s^{\alpha}g(s) = d_s^{\alpha}(s) 
= \frac{1}{\Gamma(1-\alpha)}\int^s_0 (s-\xi)^{-\alpha}\frac{dg}{d\xi}
(\xi) d\xi.
$$
Therefore, change of the order of integration yields
\begin{align*}
S &= \frac{1}{\Gamma(1-\alpha)}
\int^t_0\MLML \left( \int^s_0 (s-\xi)^{-\alpha}\frac{dg}{d\xi}
(\xi) d\xi\right)\, ds\\
=&  \frac{1}{\Gamma(1-\alpha)} \int^t_0
\frac{dg}{d\xi}(\xi) \left( \int^t_{\xi}
\MLML (s-\xi)^{-\alpha}\, ds \right) d\xi\\
=& \int^t_0 \frac{dg}{d\xi}(\xi) 
\frac{1}{\Gamma(1-\alpha)} \left(
\int^{t-\xi}_0 \eta^{\alpha-1}\MLAA(-\la_n\eta^{\alpha})
(t-\xi-\eta)^{-\alpha} d\eta \right) d\xi.
\end{align*}
For the last equality we used the change of variables $s \to \eta$ by 
$\eta = t-s$.
Moreover,
$$
\frac{1}{\Gamma(1-\alpha)} 
\int^{t-\xi}_0 \eta^{\alpha-1}\MLAA(-\la_n\eta^{\alpha})
(t-\xi-\eta)^{-\alpha} d\eta = \MLO(-\la_n(t-\xi)^{\alpha})
$$
(e.g., formula (1.100) (p.25) in \cite{Po}).
Hence, again applying integration by parts, we obtain
\begin{align*}
 S  &=\int^t_0 \frac{dg}{d\xi}(\xi)\MLO(-\la_n (t-\xi)^{\alpha}) d\xi \\
&= \Bigl[g(\xi)\MLO(-\la_n (t-\xi)^{\alpha})\Bigr]^{\xi=t}_{\xi=0}
- \int^t_0 g(\xi) \frac{d}{d\xi} \MLO(-\la_n (t-\xi)^{\alpha}) d\xi.
\end{align*}
Now, by the definition of the Mittag-Lefller function in view of 
the power series, the termwise differentiation yields
\begin{equation}\label{eqn:diff_ML} %\eqno{(2.5)} but should be (2.6)
\frac{d}{d\xi} \MLO(-\la_n (t-\xi)^{\alpha})
= \la_n(t-\xi)^{\alpha-1} \MLAA (-\la_n (t-\xi)^{\alpha}),
\quad 0<\xi<t<T.             
\end{equation}
Therefore, using $g(0) = 0$ by $g \in \CCO$, we have
$$ 
S = g(t) - \int^t_0 \la_n \MLML g(s)\, ds.
$$
Substituting this into the above we obtain
\begin{align*}
& \int^t_0 \MLML \langle f(\cdot,s), \va_n\rangle s\, ds
= -g(t)\bigl\langle \frac{x^2}{2}, \va_n\bigr\rangle\\
+& \int^t_0 \MLML g(s)\, ds
\Bigl( \bigl\langle \la_n\frac{x^2}{2}, \,\va_n\bigr\rangle
+ \bigl\langle 1+\frac{x^2}{2}p, \,\va_n\bigr\rangle\Bigr).
\end{align*}

Here by integration by parts, we calculate
\begin{align*}
& \bigl\langle \la_n\frac{x^2}{2}, \,\va_n\bigr\rangle
+ \bigl\langle \frac{x^2}{2}p, \,\va_n\bigr\rangle
= \bigl\langle \la_n\va_n +p\va_n, \frac{x^2}{2}\bigr\rangle \\
&= \bigl\langle -\frac{d^2\va_n}{dx^2}, \, \frac{x^2}{2}\bigr\rangle
= \Bigl[ -\frac{d\va_n}{dx}(x)\frac{x^2}{2}\Bigr]^{x=1}_{x=0}
+ \int^1_0 x\frac{d\va_n}{dx}(x) dx\\
&= \Bigl[x\va_n(x)\Bigr]^{x=1}_{x=0} - \int^1_0 \va_n(x)\, dx
= 1 - \langle \va_n, 1\rangle.
\end{align*}

Hence
$$
\bigl\langle \la_n\frac{x^2}{2}, \,\va_n\bigr\rangle
+ \bigl\langle 1+\frac{x^2}{2}p, \,\va_n\bigr\rangle
= 1 - \langle\va_n,1\rangle + \langle 1, \va_n\rangle = 1,
$$
so that 
\begin{align*}
& \int^t_0 \MLML \langle f(\cdot,s), \va_n\rangle\, ds\\
&\qquad = -\bigl\langle \frac{x^2}{2}, \,\va_n\big\rangle g(t)
+ \int^t_0 \MLML g(s)\, ds.
\end{align*}
Substituting this into \eqref{eqn:v-representation},
 since $\left\{ \frac{1}{\sqrt{\rho_n}}\va_n \right\}_{n\in \N}$ 
is an orthonormal basis in $L^2(0,1)$, we see
\begin{align*}
& v(x,t) = -\sumn \frac{1}{\rho_n}\bigl\langle \frac{x^2}{2},
\,\va_n\bigr\rangle g(t) \va_n(x)
+ \sumn \frac{1}{\rho_n}\int^t_0 \MLML g(s)\, ds \va_n(x)\\
&= -\frac{x^2}{2}g(t) + \sumn \frac{1}{\rho_n}\left(\int^t_0 \MLML g(s)\, ds
\right)\,\va_n(x).
\end{align*}
Since $u=v+\frac{x^2}{2}g(t)$, we have proved
\eqref{eqn:u_p,alpha-representation} for $g \in \CCO$.

Next we have to prove \eqref{eqn:u_p,alpha-representation} for
$g \in H_{\alpha}(0,T)$.
In equations \eqref{eqn:v-f} and \eqref{eqn:v-representation}, we write
$f := f_g$ and $v := v_g$ respectively in order to specify the dependence on $g$.
Since $\CCO := \{ h\in C^1[0,T];\, h(0) = 0\}$ is dense in 
$H_{\alpha}(0,T)$ (e.g., Lemma 2.2 in \cite{KRY}), 
for each $g \in H_{\alpha}(0,T)$, we can find a sequence
$g_{\ell} \in \CCO$, $\ell \in \N$ such that $g_{\ell} \longrightarrow
g$ in $H_{\alpha}(0,T)$.
Then, since $\pppa g_{\ell} \longrightarrow \pppa g$ in $L^2(0,T)$
(e.g., Theorem 2.4 in \cite{KRY}), it follows that 
$f_{g_{\ell}} \longrightarrow f_g$ in $L^2(0,T)$.
Therefore, applying the well-posedness for the initial boundary 
value problem (e.g., \cite{GLY}, Theorem 4.1 in \cite{KRY},
\cite{SY}), we see that $v_{g_{\ell}} \longrightarrow v_g$ in 
$H_{\alpha}(0,T;L^2(0,1)) \cap L^2(0,T;H^2(0,1))$.

As we already proved for $g_{\ell} \in \CCO$, we have
\begin{equation}\label{eqn:cgt_to_vg}  %\eqno{(2.6)}
\sumn \frac{1}{\rho_n} \int^t_0 \MLML g_{\ell}(s)\, ds \va_n\,
- \frac{x^2}{2}g_{\ell}(t) \longrightarrow v_g          
\end{equation}
in the space $\;H_{\alpha}(0,T;L^2(0,1)) \cap L^2(0,T;H^2(0,1))$.

On the other hand, let $h \in L^2(0,T)$.
Then one can prove by the asymptotic behavior of $\va_n$ for large 
$n \in \N$ (e.g., Section 2 of Chapter 1 of Levitan and 
Sargsjan \cite{LS}), 
that there exists a constant $\rho_0 > 0$ such that 
\begin{equation}\label{eqn:rhon>rho0}  %\eqno{(2.7)}
\rho_n \ge \rho_0 \quad \mbox{for all $n \in \N$}. 
\end{equation}
Henceforth $C>0$ denotes generic constants which are independent of 
$n$ and  choices of $h, g, t\in (0,T)$.
Let $\psi\in C^{\infty}_0 ((0,1) \times (0,T))$.
Then by integration by parts
$$
\langle\va_n, \psi(\cdot,s)\rangle 
= \frac{1}{\la_n} \langle\la_n\va_n, \, \psi(\cdot,s)\rangle
= \frac{1}{\la_n}\langle A_p\va_n, \, \psi(\cdot,s)\rangle 
= \frac{1}{\la_n} \langle\va_n, A_p\psi(\cdot,s)\rangle.
$$ 
Therefore, 
\begin{align*}
& \Bigl\langle
\sumn \frac{1}{\rho_n} \left(\int^t_0 \MLML h(s)\, ds\right) \va_n, \,\,
\psi\Bigr\rangle_{L^2((0,1)\times (0,T))}\\
&\quad = \sumn \frac{1}{\rho_n} \Bigl\langle \int^t_0 \MLML h(s)\, ds, \,\,
\frac{1}{\la_n} \langle\va_n, \, A_p\psi(\cdot,t)\rangle_{L^2(0,1)} 
\Bigr\rangle_{L^2(0,T)}.
\end{align*}
Hence, also by \eqref{eqn:rhon>rho0} and the Cauchy-Schwarz inequality, 
we have
\begin{align*}
& \left\vert \Bigl\langle
\sumn \frac{1}{\rho_n} \left(\int^t_0 \MLML h(s)\, ds \right)\va_n, \,\,
\psi\Bigr\rangle_{L^2((0,1)\times (0,T))}\right\vert\\
&\quad\le C\sumn \frac{1}{\la_n} 
\left\Vert \int^t_0 \MLML h(s)\, ds \right\Vert_{L^2(0,T)}
\Vert A_p\psi\Vert_{L^2(0,T;L^2(0,1))} \frac{\Vert \va_n\Vert}{\rho_n}\\
\le& C\sumn \frac{1}{\la_n} 
\Vert s^{\alpha-1}\MLAA(-\la_ns^{\alpha}) * h\Vert_{L^2(0,T)}.
\end{align*}
Here and henceforth we set $(g_1*g_2)(t) := \int^t_0 g_1(t-s)g_2(s)\, ds$.
By a bound of $\MLAA(-\la_ns^{\alpha})$ (e.g., Theorem 1.6 (p.35) in 
\cite{Po}), we have $\vert \MLAA(-\la_ns^{\alpha})\vert \le C$ for 
all $n\in \N$ and $s>0$.  Hence, Young's inequality yields
$$
\Vert s^{\alpha-1}\MLAA(-\la_ns^{\alpha}) * h\Vert_{L^2(0,T)}
\le \Vert s^{\alpha-1}\MLAA(-\la_ns^{\alpha})\Vert_{L^1(0,T)}
\Vert h\Vert_{L^2(0,T)}
\le C\Vert h\Vert_{L^2(0,T)}.
$$
Since $C^{-1}n^2 \le \la_n \le Cn^2$ for all $n \in \N$ (e.g., \cite{LS}),
we can obtain
\begin{align*}
&\left\vert \Bigl\langle
\sumn \frac{1}{\rho_n} \left(\int^t_0 \MLML h(s) \,ds \right)\va_n, \,\,
\psi\Bigr\rangle_{L^2((0,1)\times (0,T))}\right\vert\\
&\qquad \le C\sumn \frac{1}{n^2} \Vert h\Vert_{L^2(0,T)}
\le C\Vert h\Vert_{L^2(0,T)}
\end{align*}
for all $\psi \in C^{\infty}_0((0,1) \times (0,T))$.

Therefore, setting $h:= g- g_{\ell}$, we see that 
\begin{align*}
& \sumn \frac{1}{\rho_n} \int^t_0 \MLML g_{\ell}(s)\, ds\, \va_n\,
- \frac{x^2}{2}g_{\ell}(t)\\
&\quad\longrightarrow\,
\sumn \frac{1}{\rho_n} \int^t_0 \MLML g(s)\, ds\, \va_n\,
- \frac{x^2}{2}g(t) \quad \mbox{in $(C_0^{\infty}((0,1) \times (0,T)))'$.}
\end{align*}
In view of \eqref{eqn:cgt_to_vg}, the convergence is in 
$H_{\alpha}(0,T;L^2(0,1)) \cap L^2(0,T;H^2(0,1))$, and 
both limits in \eqref{eqn:cgt_to_vg} and the above must coincide.  Hence,
$$
v_g(x,t) = \sumn \frac{1}{\rho_n} \int^t_0 \MLML g(s)\,ds\, \va_n\,
- \frac{x^2}{2}g(t)
$$
in $H_{\alpha}(0,T;L^2(0,1)) \cap L^2(0,T;H^2(0,1))$.
Since $u_{p,\alpha}(x,t) = v_g(x,t) + \frac{x^2}{2}g(t)$ by \eqref{eqn:v-f},
the proof of Proposition 3 is complete.
$\square$
\\

We conclude this section with the following lemma.
\\
{\bf Lemma 4.}\\
{\it
Let $K_{p,\alpha}(x,t)$ be defined by
$$
K_{p,\alpha}(x,t) 
:= \sumn \frac{\va_n(x)}{\rho_n} 
\int^t_0 s^{\alpha-1}\MLAA(-\la_ns^{\alpha})\, ds
= \sumn \frac{\va_n(x)}{\la_n\rho_n} (1 - \MLO(-\la_n t^{\alpha}))
$$
$$
\mbox{for all $x\in [0,1]$ and $t\in [0,T]$}.
$$
Then,
\\
(i) The series is uniform convergent in $x \in [0,1]$ and $t \in [0,T]$,
and $K_{p,\alpha}(x,\cdot) \in L^{\infty}(0,\infty)$ and is analytic in $t>0$
for all $x\in [0,1]$.
\\
(ii)
$$
\int^{\xi}_0 u_{p,\alpha}(x,t)\, dt = (K_{p,\alpha}(x,\cdot) * g)(\xi) \quad
\mbox{for all $x\in [0,1]$ and $\xi\in [0,T]$}.
$$
}
\\
{\bf Proof of (i).}\\
In view of \eqref{eqn:diff_ML}, we have
\begin{equation}\label{eqn:mlf-estimates}  %\eqno{(2.8)}
\int^t_0 s^{\alpha-1}\MLAA(-\la_ns^{\alpha})\,ds
= \frac{1}{\la_n} \int^0_t \frac{d}{ds}(\MLO(-\la_ns^{\alpha}))\, ds
= \frac{1}{\la_n}(1 - \MLO(-\la_nt^{\alpha})).
\end{equation}
Hence, 
$$
K_{p,\alpha}(x,t) = \sumn \frac{1}{\la_n\rho_n} (1 - \MLO(-\la_n t^{\alpha}))
\va_n(x),  \quad 0\le x \le 1, \, t>0.
$$
From Theorem 1.6 (p.35) in \cite{Po}, we know that there exist
constants $C>0$ and $\theta_0>0$ such that 
$$
\vert \MLO(-\la_nz^{\alpha})\vert \le C \quad 
\mbox{for all $n\in \N$ and $z \in \Sigma := \{ z\in \C;\, \,
\vert \mbox{Arg}\,z \vert < \theta_0\}$}.
$$
We fix small $\delta>0$ arbitrarily.  Since 
$\Vert \va_n\Vert_{H^{\theta}(0,1)} \le C\Vert A_p^{\frac{\theta}{2}}\va_n
\Vert_{L^2(0,1)}$ with $0<\theta<2$, applying the Sobolev embedding and
recalling $\rho_n = \Vert\va_n\Vert^2_{L^2(0,1)}$, we have
$$
\Vert \va_n\Vert_{C[0,1]} \le C\Vert \va_n\Vert_{H^{\frac{1}{2}+\delta}(0,1)}
\le C\Vert A_p^{\frac{1}{4}+\frac{\delta}{2}}\va_n\Vert_{L^2(0,1)}
= C\la_n^{\frac{1}{4}+\frac{\delta}{2}}\sqrt{\rho_n}.
$$
Hence, by (2.8), we obtain
$$
\left\vert \frac{1}{\rho_n\la_n} (1 - \MLO(-\la_nz^{\alpha}))\va_n(x)
\right\vert 
\le \frac{C}{\la_n\sqrt{\rho_n}}\la_n^{\frac{1}{4}+\frac{\delta}{2}},
\quad 0\le x\le 1, \, \, z\in \Sigma, 
$$
and so
$$
\sumn \frac{1}{\la_n\rho_n} \vert (1 - \MLO(-\la_nz^{\alpha}))
\va_n(x)\vert
$$
$$
\le C\sumn \frac{1}{\la_n^{\frac{3}{4}-\frac{\delta}{2}}}
\le C\sumn \frac{1}{n^{\frac{3}{2}- \delta}} < \infty,
\quad 0\le x\le 1, \, \, z\in \Sigma.         \eqno{(2.10)}
$$
Here we used $\la_n \sim n^2$ (e.g., \cite{LS}).
Since $\MLO(-\la_nz^{\alpha})$ is analytic in $z \in \Sigma$, we can complete
the proof of (i).
\\
{\bf Proof of (ii).}\\
Since the series in (2.3) is convergent in $L^2(0,T;H^2(0,1))$, 
by $H^2(0,1) \subset C[0,1]$, we see that 
$$
u_{p,\alpha}(x,t) = \int^t_0 \left( \sumn \frac{1}{\rho_n} 
\MLML \va_n(x) \right) g(s) \,ds
$$
is convergent in $L^2(0,T;C[0,T])$.  Therefore,
$$
\int^{\xi}_0 u_{p,\alpha}(x,t)\, dt 
= \int^{\xi}_0 \left\{\int^t_0 \left( \sumn \frac{1}{\rho_n}
\MLML \va_n(x)\right) g(s) \,ds \right\}\, dt
$$
for all fixed $x\in [0,1]$.
Exchanging the orders of the integrals and changing the variables
$t \to \eta$: $\eta = t-s$, we obtain
\begin{align*}
& \int^{\xi}_0 \left( \int^t_0
\MLML g(s)\, ds \right)\, dt\\
=& \int^{\xi}_0 \left( \int^{\xi}_s
\MLML\, dt \right) g(s)\, ds
= \int^{\xi}_0 \left( \int^{\xi-s}_0
\eta^{\alpha-1}\MLAA(-\la_n\eta^{\alpha})\, d\eta \right) g(s)\, ds.
\end{align*}
Hence by (2.9), we have verified (ii)  and the proof of Lemma~4 is 
complete.
\section{Proof of Theorem 2}

Let 
$$
\left\{ \begin{array}{rl}
&A_qw(x) = - \frac{d^2w}{dx^2}(x) - q(x)w(x), \quad 0<x<1,\\
& \mathcal{D}(A_q) = \left\{ w \in H^2(0,1);\,
\frac{dw}{dx}(0) = \frac{dw}{dx}(1) = 0 \right\}.
\end{array}\right.
$$
We let
$ 0 < \mu_1 < \mu_2 < \cdots,$
denote all the eigenvalues of the operator $A_q$ and let
$\psi_n$, $n\in \N$ be the corresponding eigenfunction for $\mu_n$,
that is  $\psi_n \in \mathcal{D}(A_q)$ satisfies
$A_q\psi_n = \mu_n\psi_n$ in $(0,1)$ and we take the normalization
of the eigenfunctions to be $\psi_n(1) = 1$.
Given this, we set $\, \sigma_n := \Vert \psi_n\Vert^2$, for $n\in \N$.
Similarly to Lemma 4, we define
$$
K_{q,\beta}(1,t) = \sumn \frac{1}{\sigma_n}
\int^t_0 s^{\beta-1}E_{\beta,\beta}(-\mu_ns^{\beta})\,ds, \quad 
t>0.
$$
It is sufficient to prove the theorem with data 
$u_{p,\alpha}(1,t) = u_{q,\beta}(1,t)$, $0<t<T$.  For the other case
at $x=0$, replacing the conditions $\va_n(1) = \psi_n(1) = 1$ by 
$\va_n(0) = \psi_n(0) = 1$, we can repeat the whole argument and thus
omit the details for this case.

Since $\int^{\xi}_0 u_{p,\alpha}(1,t)\, dt
= \int^{\xi}_0 u_{q,\beta}(1,t)\, dt$ by $u_{p,\alpha}(1,t) 
= u_{q,\beta}(1,t)$ for $0<t<T$, in view of Lemma 4, we see
$$
(K_{p,\alpha}(1,\cdot) * g)(\xi) = (K_{q,\beta}(1,\cdot) * g)(\xi), 
\quad 0<\xi<T,
$$
that is, 
$$
((K_{p,\alpha} - K_{q,\beta})(1,\cdot) * g)(t)=0, \quad 0<t<T.
$$
Since $g\not \equiv 0$, we apply the Titchmarsh convolution theorem (e.g.,
Titchmarsh \cite{Ti}), so that there exists $t_0>0$ such that 
$$
K_{p,\alpha}(1,t) = K_{q,\beta}(1,t), \quad 0<t<t_0.
$$
Lemma 4 implies that $K_{p,\alpha}(1,t)$ and $K_{q,\beta}(1,t)$ 
are analytic in $t>0$, and so
$$
\sumn \frac{1}{\la_n\rho_n} (1 - \MLO(-\la_nt^{\alpha}))
= \sumn \frac{1}{\mu_n\sigma_n} (1 - E_{\beta,1}(-\mu_nt^{\beta})),
\quad t>0.
$$
By the asymptotics of $\MLO(-\eta)$ and $E_{\beta,1}(-\eta)$ for 
large $\eta>0$ (e.g., Theorem 1.4 (pp.33-34) in \cite{Po}), we have
$$
\MLO(-\la_nt^{\alpha}) = \frac{1}{\Gamma(1-\alpha)}\frac{1}{\la_nt^{\alpha}}
+ O\left(\frac{1}{t^{2\alpha}}\right)
$$
and
$$
E_{\beta,1}(-\mu_nt^{\beta}) = \frac{1}{\Gamma(1-\beta)}\frac{1}
{\mu_nt^{\beta}}
+ O\left(\frac{1}{t^{2\beta}}\right)
$$
for all large $t>0$.
Hence
\begin{align*}
& \sumn \frac{1}{\la_n\rho_n} 
- \frac{1}{\Gamma(1-\alpha)}\sumn \frac{1}{\la_n\rho_n}
\frac{1}{\la_nt^{\alpha}} + O\left(\frac{1}{t^{2\alpha}}\right)\\
=& \sumn \frac{1}{\mu_n\sigma_n} 
- \frac{1}{\Gamma(1-\beta)}\sumn \frac{1}{\mu_n\sigma_n}
\frac{1}{\mu_nt^{\beta}} + O\left(\frac{1}{t^{2\beta}}\right)
\end{align*}
for large $t>0$.  Letting $t \to \infty$, we obtain 
$$
\sumn \frac{1}{\la_n\rho_n} = \sumn \frac{1}{\mu_n\sigma_n}.
$$
Assume that $\alpha > \beta$.  Then
$$
- \frac{1}{\Gamma(1-\alpha)}\sumn \frac{1}{\la_n\rho_n}
\frac{1}{\la_nt^{\alpha-\beta}} + O\left(\frac{1}{t^{2\alpha-\beta}}\right)
= - \frac{1}{\Gamma(1-\beta)}\sumn \frac{1}{\mu_n\sigma_n}
\frac{1}{\mu_n} + O\left(\frac{1}{t^{\beta}}\right)
$$
for large $t>0$. 
Letting $t \to \infty$, we obtain
$$
\frac{1}{\Gamma(1-\beta)}\sumn \frac{1}{\mu_n^2\sigma_n} = 0.
$$
Since $\sigma_n = \Vert \psi_n\Vert^2 > 0$, this is impossible.
Hence $\alpha \le \beta$.
By an entirely similar argument we see that $\alpha < \beta$ is impossible
and so conclude that  $\alpha = \beta$.

Now we move to complete the proof of the theorem. 
We see
\begin{equation}\label{eqn:mlf_lambda_mu}  %\eqno{(3.1)}
\sumn \frac{1}{\la_n\rho_n} \MLO(-\la_n t^{\alpha})
= \sumn \frac{1}{\mu_n\sigma_n} E_{\alpha,1}(-\mu_n t^{\alpha}),
\quad t>0.    
\end{equation}
Now we can argue similarly to \cite{CNYY}.
Using 
$$
\left\vert \frac{1}{\la_n\rho_n} \MLO(-\la_nt^{\alpha})\right\vert
\le \frac{C}{\la_n}, \quad n\in \N, \quad t>0,
$$
we see that the series in (3.1) are convergent uniformly in $[0, \infty)$.
Therefore we can take the Laplace transforms termwise to have
$$
\sumn \frac{1}{\la_n\rho_n} \int^{\infty}_0 e^{-\zeta t}
\MLO(-\la_n t^{\alpha}) \,dt
= \sumn \frac{1}{\mu_n\sigma_n} \int^{\infty}_0 e^{-\zeta t}
E_{\alpha,1}(-\mu_n t^{\alpha}) \,dt, \quad \zeta>0.  
$$
By formula (1.80) (p.21) in \cite{Po}, we obtain
$$
\sumn \frac{1}{\la_n\rho_n} \frac{\zeta^{\alpha-1}}{\zeta^{\alpha}
+ \la_n}
= \sumn \frac{1}{\mu_n\sigma_n} \frac{\zeta^{\alpha-1}}{\zeta^{\alpha}
+ \mu_n}, \quad \zeta>0.
$$
Dividing by $\zeta^{\alpha-1}$ and setting $\eta = \zeta^{\alpha}$, we have
\begin{equation}\label{eqn:mlf_lambda_mu2}  %\eqno{(3.2)}
\sumn \frac{1}{\la_n\rho_n} \frac{1}{\eta + \la_n}
= \sumn \frac{1}{\mu_n\sigma_n} \frac{1}{\eta + \mu_n}, \quad \eta>0.
\end{equation}
Since $\la_n \sim n^2$ and $\mu_n \sim n^2$ for large $n\in \N$, we see 
that both sides of (3.2) are convergent uniformly in any compact set 
in $\C \setminus (\{-\la_n\}_{n\in \N} \cup \{-\mu_n\}_{n\in \N})$ and 
are analytic in $\C \setminus (\{-\la_n\}_{n\in \N} \cup 
\{-\mu_n\}_{n\in \N})$.

Assume that $\la_m \not\in \{\mu_n\}_{n\in \N}$ for $m\in\N$.  
Then we can choose a small circle $C_m$ centered at $-\la_m$ and 
$\{-\mu_n\}_{n\in \N}$ is not included in the disk centered at $-\la_n$ 
bounded by $C_m$.  Integrating on $C_m$ and applying the 
Cauchy theorem, we have
$$
\frac{2\pi \sqrt{-1}}{\la_m\rho_m} = 0,
$$
which is impossible. Hence $\la_m \in \{\mu_n\}_{n\in \N}$ for each
$m$.  Similarly $\mu_m \in \{\la_n\}_{n\in \N}$ for each $m\in \N$.
Therefore
\begin{equation}\label{eqn:lambda=mu}  %\eqno{(3.3)}
\la_n = \mu_n, \quad n\in \N. 
\end{equation}
By \eqref{eqn:mlf_lambda_mu2}, we have
$$
\sumn \left( \frac{1}{\la_n\rho_n} - \frac{1}{\la_n\sigma_n}
\right) \frac{1}{\eta + \la_n} = 0, \quad
\eta \in \C \setminus \{-\la_n\}_{n\in \N}.
$$
Again integrating on $C_m$, we obtain
$$
\frac{2\pi\sqrt{-1}}{\la_n}\left( \frac{1}{\rho_n} - \frac{1}{\sigma_n}\right)
= 0,
$$
that is,
\begin{equation}\label{eqn:rho=sigma}  %\eqno{(3.4)}
\rho_n = \sigma_n, \quad n\in \N.   
\end{equation}
Now, using \eqref{eqn:lambda=mu} and \eqref{eqn:rho=sigma},
we apply the Gel'fand-Levitan theory (e.g., Gel'fand and Levitan
\cite{GL}), and we can obtain $p(x) = q(x)$ for $0<x<1$.  
The application is similar to \cite{CNYY}, \cite{Mu}, \cite{SM}, and so 
we omit the details.
Thus the proof of Theorem~2 is complete.
$\square$

\section{Application of the representation formula}

The representation formula Proposition 3 is useful for qualitative analyses
of fractional equations. 
Here we explain one application.
 
We let $0<\alpha<1$ and we fix $p\in C[0,1]$, $\le 0$ on $[0,1]$.
Let
$$
\left\{ \begin{array}{rl}
& \pppa u(x,t) = \ppp_x^2u(x,t) + p(x)u(x,t), \quad 0<x<1,\, 0<t<T,\\
& u(0,t) = \ppp_xu(0,t) = 0, \quad 0<t<T, \\
& u \in H_{\alpha}(0,T;L^2(0,1)).
\end{array}\right.
                                         \eqno{(4.1)}
$$
Then we are interested in the question: can we conclude $u(x,t) = 0$
for $0<x<1$ and $0<t<T$?
\\

This is a kind of unique continuation property under the assumption 
$u \in H_{\alpha}(0,T;L^2(0,1))$ which can be interpreted as that an 
initial value of $u$ is zero. This kind of unique continuation was proved 
by Cheng, Lin and Nakamura \cite{CLN} for $\alpha = \frac{1}{2}$,
Lin and Nakamura \cite{LN1} for $\alpha\in (0,1)$ and
Lin and Nakamura \cite{LN2} for $\alpha\in (0,1) \cup (1,2)$ for 
general time-fractional partial differential equations.
Their proofs are based on the techniques of pseudo-differential 
operators.
 
For $\alpha=1$, we can prove the unique continuation without any information 
of initial conditions, and the corresponding unique continuation is
proved for a one-dimensional time-fractional equation by Li and Yamamoto 
\cite{LiYa2019}.  More precisely, if $u$ is in a suitable class and 
satisfies  
$$
\left\{ \begin{array}{rl}
& \pppa u(x,t) = \ppp_x^2u(x,t), \quad 0<x<1,\, 0<t<T,\\
& u(0,t) = \ppp_xu(0,t) = 0, \quad 0<t<T, 
\end{array}\right.
$$
then $u(x,t) = 0$ for $0<x<1$ and $0<t<T$. 

However such unique continuation not requiring any initial conditions,
is not known for general case in multidimensions.

In this section, for the one-dimensional case (4.1),
we provide a simpler proof than \cite{CLN, LN1, LN2}, 
which relies on the representation formula Proposition 3. 
\\
{\bf Proposition 5.}\\
{\it 
Let $u \in H_{\alpha}(0,T;L^2(0,1)) \cap L^2(0,T;H^2(0,1))$ satisfy 
(4.1) and $\ppp_xu(1,\cdot) \in H_{\alpha}(0,T)$.
Then $u(x,t) = 0$, $0<x<1$, $0<t<T$.
}
\\

By the definition of $H_{\alpha}(0,T)$ given in Section 1, 
if $0<\alpha < \frac{1}{2}$, then $H_{\alpha}(0,T) 
= H^{\alpha}(0,T)$ and in (4.1) the condition $u \in H_{\alpha}(0,T;L^2(0,1))$
does not require anything for the behavior of the solution $u$ near $t=0$.
In other words, we need not pose any conditions at $t=0$ to $u$.

It seems that we can remove a condition $\ppp_xu(1,\cdot) \in H_{\alpha}(0,T)$,
but we here omit the details.
\\
\vspace{0.2cm}
\\
{\bf Proof.}\\
We set $g:= \ppp_xu(1,\cdot) \in H_{\alpha}(0,T)$.
Then $u$ satisfies (1.4) and $u(0,t) = 0$, $0<t<T$.
Lemma 4 (ii) implies
$$
(K_{p,\alpha}(0,\cdot) * g)(t) = 0, \quad 0<t<T.
$$
By the Titchmarsh theorem on the convolution (e.g., \cite{Ti}), there
exist $t_1, t_2 \ge 0$ such that
$$\left\{ \begin{array}{rl}
& t_1 + t_2 = T, \\
& K_{p,\alpha}(0,t) = 0, \quad 0\le t \le t_1, \\
& g(t) = 0, \quad 0\le t \le t_2.
\end{array}\right.
                                   \eqno{(4.2)}
$$
Assuming that $t_1 > 0$, we will derive a contradiction, which proves
$t_1=0$, that is, $g=0$ in $(0,T)$.  The argument is similar to 
the proof of Theorem 2.

The analyticity of $K_{p,\alpha}(0,t)$ in 
$t>0$ yields $K_{p,\alpha}(0,t) = 0$ for all $t>0$.

Since $\lim_{t\to \infty} \MLO(-\la_nt^{\alpha}) = 0$ 
(e.g., Theorem 1.6 (p.35) in \cite{Po}), we have 
$$
\lim_{t\to \infty} K_{p,\alpha}(0,t) = \sumn \frac{\va_n(0)}{\la_n\rho_n}
= 0.
$$
Hence
$$
\sumn \frac{\va_n(0)}{\la_n\rho_n}\MLO(-\la_nt^{\alpha}) = 0, \quad
t>0.
$$
This series is convergent in $L^{\infty}(0,\infty)$ and so we can take
the Laplace transform term by term.  In view of formula (1.80) (p.21) in 
\cite{Po}, we obtain
$$
\sumn \frac{\va_n(0)}{\la_n\rho_n}\frac{z^{\alpha-1}}{z^{\alpha}+\la_n}
= 0, \quad \mbox{Re}\, z > 0.
$$
Then, dividing by $z^{\alpha-1}$ and setting $\eta = z^{\alpha}$, we 
have
$$
\sumn \frac{\va_n(0)}{\la_n\rho_n}\frac{1}{\eta+\la_n} = 0, \quad 
\mbox{Re}\, \eta > 0.
$$
Similarly to (2.10), we can verify
$$
\sumn \frac{\va_n(0)}{\la_n\rho_n} < \infty,
$$
and so we can continue analytically in $\eta$ as much as possible to 
obtain
$$
\sumn \frac{\va_n(0)}{\la_n\rho_n}\frac{1}{\eta+\la_n} = 0, \quad 
\eta \in \C \setminus \{-\la_n\}_{n\in \N}.
$$
Choosing a small circle $\Gamma_1$ centered at $-\la_1$ such that 
the interior of the disk bounded by $\Gamma_1$ does not contain 
$-\la_n$ with $n\ge 2$ and integrating on $\Gamma_1$, in terms of 
the Cauchy theorem, we see
$$
\frac{\va_1(0)}{\la_1\rho_1}2\pi \sqrt{-1} = 0,
$$
that is, $\va_1(0) = 0$.  Since $\frac{d^2\va_1}{dx^2}(x) 
+ (p(x)+\la_1)\va_1(x) = 0$, $0<x<1$ and $\frac{d\va_1}{dx}(0) = 0$,
we have $\va_1(x) = 0$ for all $0<x<1$, which is impossible.

Then we can conclude that $t_1=0$.  By (4.2), we reach $g(t)=0$ for
$0<t<T$.  Thus the proof of Proposition 5 is complete.
 
\section*{Acknowledgments}

The first author  was supported in part by the
National Science Foundation through award {\sc dms}-1620138.

The second author was supported by Grant-in-Aid for Scientific Research (S)
15H05740 of Japan Society for the Promotion of Science and
by The National Natural Science Foundation of China
(no. 11771270, 91730303).
This work was prepared with the support of the "RUDN University Program 5-100".

%%%%%%%%%%%%%%%%%%%%%%%%%%%%%%%%%%%%%%

\end{document}